\documentclass[final]{siamltex}

\pdfoutput = 1

\usepackage{graphicx}
\usepackage{psfrag}
\usepackage{amsmath,rotating}
\usepackage{amstext,amssymb,epsf}
\usepackage[dvips]{epsfig}
\usepackage{subfigure}
\usepackage[numbers]{natbib}
\usepackage{placeins}
\usepackage{amsbsy} 

\newcommand{\St}{\text{St}}
\newcommand{\taup}{\tau_p}

\newcommand{\tf}{\overline{f}}
\newcommand{\xx}{\mathbf{x}}
\newcommand{\vv}{\mathbf{v}}

\newcommand{\uu}{\mathbf{u}}

\newcommand{\FF}{\mathbf{F}}

\title{On the Eulerian Large Eddy Simulation of disperse phase flows:  an asymptotic preserving scheme for small Stokes number flows}

\author{C.~Chalons\footnotemark[1] \footnotemark[2], M.~Massot\footnotemark[2] \footnotemark[3] \footnotemark[4], \and A.~Vi\'{e}\footnotemark[3] \footnotemark[4] \footnotemark[5]}

\begin{document}

\maketitle
{\begin{center} \today \end{center}}
\renewcommand{\thefootnote}{\fnsymbol{footnote}}
\footnotetext[1]{Laboratoire de Math\'ematiques de Versailles, UMR~8100,
Universit\'e de Versailles Saint-Quentin-en-Yvelines,
UFR des Sciences, b\^atiment Fermat, 45 avenue des Etats-Unis, 78035 Versailles cedex}
\footnotetext[2]{F\'ed\'eration de Math\'ematiques - FR CNRS 3487, Ecole Centrale Paris, Grande Voie des Vignes, 92295 Chatenay-Malabry, France }
\footnotetext[3]{CNRS Laboratoire EM2C - UPR 288, Grande Voie des Vignes, 92295 Chatenay-Malabry, France}
\footnotetext[4]{Ecole Centrale Paris, Grande Voie des Vignes, 92295 Chatenay-Malabry, France }
\footnotetext[5]{Center For Turbulence Research, Stanford University, 488, Escondido Mall, Stanford CA 94305-3035, USA}
\renewcommand{\thefootnote}{\arabic{footnote}}

\begin{abstract}

In the present work, the Eulerian Large Eddy Simulation of dilute disperse phase flows is investigated. By highlighting the main advantages and drawbacks of the available approaches in the literature, a choice is made in terms of modelling:  a Fokker-Planck-like filtered kinetic equation proposed by Zaichik et al. 2009 and a Kinetic-Based Moment Method (KBMM) based on a Gaussian closure for the NDF proposed by Vie et al. 2014. The resulting Euler-like system of equations is able to reproduce the dynamics of particles for small to moderate Stokes number flows,
given a LES model for the gaseous phase, and is representative of the generic difficulties of such models. Indeed, it encounters strong constraints in terms of numerics in the small Stokes number limit, which can lead to a degeneracy of the accuracy of standard numerical methods. These constraints are: 1/as the resulting sound speed is inversely proportional to the Stokes number, it is highly CFL-constraining, and 2/the system tends to an advection-diffusion limit equation on the number density that has to be properly approximated by the designed scheme used for the whole range of Stokes numbers. Then, the present work proposes a numerical scheme that is able to handle both. Relying on the ideas introduced in a different context by Chalons et al. 2013: a Lagrange-Projection, a relaxation formulation and a HLLC scheme with source terms, we extend the approach to a singular flux as well as properly handle the energy equation. The final scheme is proven to be Asymptotic-Preserving on 1D cases comparing to either converged or analytical solutions and can easily be extended to multidimensional configurations, thus setting the path for realistic applications.

\end{abstract}



\pagestyle{myheadings}
\thispagestyle{plain}
\markboth{A. Vi\'e, C.~Chalons, \& M.~Massot}{An AP numerical scheme for the LES of disperse phase flows}

\section{Introduction}

The simulation of disperse phase flows is nowadays of great importance in several applications, such as 
 automotive engines, aeronautical combustors or fluidized beds. Actually, the modelling of such flows relies on the accurate description of both the continuous carrier phase, gaseous or liquid, and the discrete particulate phase, composed of particles or droplets. 

In the context of small particles with respect to all carrier phase flow scales, the modelling of the carrier phase could be envisioned at a Mesoscopic level \cite{fox2012}, i.e. the flow around each particle is not resolved and the coupling effects between particles and the carrier phase are modeled using Mesoscopic closures, such as the Stokes law for drag force.
At this level, to solve the statistics of the disperse phase, a Population Balance Equation (PBE) on the Number Density Function (NDF) can be used. The NDF represents the probability of having a particle at a certain position of the phase space, the phase space dimensions being the relevant properties of the particles, like their position, velocity, size, temperature...
 To solve this equation, three approaches are possible:
\begin{itemize}
\item Full resolution: the PBE is directly solved by discretizing the entire phase space. This method is the most precise for compactly supported distributions, but is too expensive for unsteady configurations 
when the phase space of the particles is too large\footnote{For example for 3D simulations where the phase space is at least 6D, 3D for the position and 3D for the velocity of the particles.};
\item Direct Monte Carlo Simulation (DMSC): the NDF is sampled by an ensemble of individual stochastic Lagrangian realizations, which are solved by means of ODEs. This approach is less expensive than the full resolution, but may also be limited by statistical convergence issues when the dimensionality of the phase space is high and many realizations are needed;
\item Moment methods: instead of solving for the NDF directly, moments of the NDF are solved, which are integrals over the phase space. By reducing the phase space to the physical space only, this method is computationally efficient. However it requires an additional effort in terms of modelling, as the integration step results in a loss of information.
\end{itemize}

In the present work, we are interested in moment methods, because of their computational efficiency with regards to other approaches. One of the main issues with moment methods is the accurate description of the velocity distribution of the particulate phase. Actually, in turbulent flows, the velocity distribution can drastically change with the inertia of the particles, which can be quantified by the Stokes number based on the Kolmogorov time scale.
For Stokes number smaller than one, the NDF is monokinetic, i.e. all particles at the same position have the same velocity, and such a distribution can be uniquely determined using zero and first order moments, i.e. density and momentum. 
For higher Stokes number, particles trajectories may cross, and the velocity distribution is no longer a unique Dirac $\delta$-function, and higher order moments are needed. To handle these higher order moments, several methods can be found in the literature, and can be split into two categories. On one hand, Algebraic-Closure-Based Moment Methods (ACBMM) \cite{alipchenkov2007,kaufmann2008,masi2014ijmf,masi2014ftc} derive closures for the second order moments using physical and/or mathematical assumptions. On the other hand, Kinetic-Based Moment Methods (KBMM) close the system by using a presumed shape for the NDF \cite{laurent2001,massot2004,ctr10b,kah10,yuan2011,vie2013cicp,vie12ctr}, which has as many parameters as the number of moments required to be controlled to describe the NDF accurately. The choice between each type of closure is motivated by the structure and the complexity of the encountered PTC, and is directly related to the number of moments.

The moment methods are a powerful tool to simulate academic configurations, but when it comes to complex configurations with a large spectrum of time and space scales, the mesh size may become too large and the computation too expensive to be achieved. To circumvent this issue, the Large Eddy Simulation is a powerful strategy: by filtering the equations in space or frequency domain, the mesh requirements can be significantly reduced. For the gas phase, the problem has been intensively studied, see \cite{smagorinsky1963,nicoud1999,nicoud2011} for example. In this work this topic is not addressed, assuming that the gas phase closures affecting the disperse phase are given, and that we are in a sufficiently dilute regime to neglect the impact of the disperse phase on the gas phase, i.e a one-way coupling regime. For the disperse phase, two types of approach can be found in the literature, which differs by the filtering procedure: 
\begin{itemize}
\item The first approach consists in deriving moment equations of the NDF, and then filtering this system to obtain the filtered moment equations. This method was investigated in \cite{moreau2010}, which proposes to use a Smagorinsky-like approach. It had already been applied in academic and industrial applications \cite{boileau08jx,martinez2009a,sanjose2011}. The main issue is that the connection between the NDF and the filtered moment equations is lost. Thus, the realizability of the moment equations, that is the moment are always those of a positive NDF,  is hard to achieved;
\item The second approach consists in filtering the kinetic equation, and then deriving moment equations on the filtered NDF. This approach was first envisioned in \cite{pandya2002}, and has been also studied in the context of the Mesoscopic Eulerian Formalism (MEF) \cite{fevrier2005} in \cite{zaichik2009}. The main interest of such an approach is that it keeps a clear link between the kinetic level, which is in fact the filtered kinetic equation, and the moment level, which is really helpful to devise realizable methods.
\end{itemize}
Whatever the method used to obtain the LES equations, one of the main issues of the LES modelling is the asymptotic behavior at small Stokes number: when the Stokes number of the particles tends to zero, the resulting system 
admits a natural  asymptotic behavior governed by an advection-diffusion 
on the number density, as the particles  become tracers diffusing in the carrier gas phase. 
This asymptotic limit is important for real applications, for example in liquid-fuel combustion systems: when the droplets evaporates, they tend to this zero Stokes number, and the right dynamics has to be captured to reproduce the right fuel distribution in the system. For now,  LES models have various way of dealing with this limit:
\begin{itemize}
\item In \cite{moreau2010}, the authors do not treat this limit, because they consider moderate Stokes number flows, for which the limit is not predominant;
\item In \cite{shotorban2007}, they obtain an advection-diffusion limit, but with an isotropic diffusion coefficient, whereas the real diffusion coefficient, which is the subgrid scale tensor of the gas phase, can be anisotropic;
\item In \cite{pandya2002} and \cite{zaichik2009}, the authors get the right limit, that is an asymptotic  advection-diffusion asymptotic limit on the density, with a potentially anisotropic diffusion coefficient.
\end{itemize}

Following these statements, in the present work we will investigate the close link between modelling and numerical methods for the LES approach of Zaichik et al. \cite{zaichik2009}. For this model, in the zero Stokes number limit, the flux and source terms that arises from the closures at the kinetic level become infinite, leading to strong constraints on the time step. Thus, using a global unified scheme for the whole Stokes range and recovering the proper asymptotic limit for small Stokes numbers is not straightforward, as these extreme constraint on the numerics will degrade the quality of the numerical approximation; standard numerical methods for hyperbolic system of conservation laws will not preserve the asymptotic limit or eventually lead to robustness issues. 


 Consequently, the goal of the present study is to propose a numerical scheme which (1) recovers the advection-diffusion asymptotic limit for the number density, referred as the asymptotic-preserving property  \cite{jin2,jinlevermore} and (2) can get rid of the small time step imposed by the source term and the acoustic waves in a standard time-explicit Godunov-type method.
To do so, we rely on the approach proposed by Chalons, Girardin and Kokh in \cite{chalons2012ap}, which is based on three ingredients:
\begin{itemize}
\item a Lagrange-Projection decomposition \cite{godlewski1996}, that
separates the terms responsible for the acoustic waves and the transport waves. This decomposition 
 overcomes the strong restriction on the time step CFL condition coming from the large value of the
sound speed when the Stokes number goes to zero by means of a time-implicit treatment of the 
Lagrangian system, see \cite{coqueletal2010} where such an approach was first proposed;
\item a relaxation strategy first introduced by Suliciu \cite{suliciu1998} and Jin and Xin \cite{jinxin}, and then further studied by many authors for instance in \cite{coqueletal2001,chalonscoquel,chalons2008,bouchut, coqueletal2010}. The motivation for using a relaxation approach is to circumvent the nonlinearities involved in 
the pressure law and then to make the time-implicit treatment of the Lagrangian step very low cost;
\item a USI (Upwinding Sources at Interfaces) approach first 
initiated by Cargo and Le Roux \cite{cargoleroux}, Greenberg and Le Roux \cite{GreenbergLeRoux}, Gosse 
and Le Roux \cite{GosseLeRoux}, 
(see also Gosse \cite{Gosse}, Perthame and Simeoni \cite{PerthameSimeoni}, 
Jin \cite{jin}, R. Botchorishvili, B.Perthame 
and A.Vasseur \cite{BotchorishviliPerthameVasseur}...)
the principle of which is to upwind the sources at interfaces in order to get the expected asymptotic-preserving 
property (or related well-balanced property). The form of the USI approach used here is the same as the one 
recently introduced in \cite{chalonsetal2010} and \cite{chalons2012ap}.
\end{itemize}
The novelty of the present contribution is related to the particular context of turbulent two-phase flows, which involve singular fluxes and infinite sound speeds in the limit as well as relaxation source terms in the energy equation, a different context compared to the original studies in the field, for which we provide a specific treatment.

The paper is organized as follow. First the modelling approach is presented, based on the kinetic equation of \cite{zaichik2009}, and on a Gaussian closure for the moment system \cite{massot2004}, arguing on the physical aspect of such a system, and presenting the asymptotic limit. Second, the main features of the numerical scheme are described: (1) a Lagrange-Projection \cite{godlewski1996} to decouple slow material and fast sound waves, which allows to use an explicit numerical method for material wave for the sake of precision and an implicit method on sound waves to avoid the CFL limitation, (2) a relaxation strategy \cite{chalons2008}, to avoid the non-linearity induced by the pressure law, and (3) an HLLC scheme which includes the source terms \cite{gallice2003}, to recover the asymptotic limit of the moment system.

\section{Moment methods for the Large Eddy Simulation of particle-laden flows}

In the present work, the following assumptions are considered:
\begin{itemize}
\item point particles: no effect of the finite size of particles. Particles are smaller than the Kolmogorow length scale of the carrier flow, so that every flow modification induced by the particles is immediately diffused by the micro-mixing of the turbulence;
\item high Knuden number: the mean free path of particles is sufficiently large to neglect collisions;
\item very dilute regime: the impact of the disperse phase on the carrier phase is neglected;
\item fixed size solid particle: no size change nor breakup are accounted for.
\end{itemize}
Consequently, in the hierarchy of modelling approaches proposed in \cite{fox2012}, here the description of the flow is done at the mesoscopic level: we consider the Number Density Function (NDF) $f(t,\xx,\vv)$ where $t$ is the time, $\xx$ the position and $\vv$ the velocity of a particle. The state of the particle is then solely described by its position and velocity.
\subsection{The Population Balance Equation and the moment methods}
The NDF $f(t,\xx,\vv)$ satisfies a Population Balance Equation:
\begin{equation}
	\partial_t f+\vv\cdot\partial_{\xx} f+\partial_{\vv}\cdot \left( \mathbf{F}_{d} f  \right)=0,
\label{WBE}
\end{equation}
where $\mathbf{F}_d$ is the drag force applied by the carrier phase on the disperse phase. The drag force is modeled using Stokes law:
\begin{equation}
\FF_d(\vv)=\dfrac{\uu_g-\vv}{\tau_p},
\end{equation}
where $\tau_p$ is the relaxation time of particles. 

In real applications, solving Eq.~(\ref{WBE}) directly is intractable for two reasons: the high dimensionality of the phase space in 3D and the wide range of scales of either carrier or disperse phases. An alternative is the macroscopic Eulerian approach: Eq.~(\ref{WBE})  is integrated over the velocity phase space, and conservation equations on moments are obtained:

\begin{align}
M_{ijk}=\iiint v_1^iv_2^jv_3^kf(t,\xx,\vv)\text{d}\vv ,
\end{align}
\begin{align}
\partial_t M_{ijk}+\partial_{\xx}\begin{pmatrix}M_{i+1,j,k} \\ M_{i,j+1,k} \\M_{i,j,k+1} \end{pmatrix}=-\dfrac{1}{\tau_p}\left((i+j+k)M_i-\uu_g\cdot \begin{pmatrix}iM_{i-1,j,k} \\ jM_{i,j-1,k} \\kM_{i,j,k-1} \end{pmatrix} \right).
\label{momenteq}
\end{align}

 This conservation equation lives now in a 3D space. However, for every moment set  of order $N=i+j+k$ that is solved, additional higher order moments of order $N+1$ are necessary. This additional moment requires a closure law, that determines (1) the accuracy of the method (2) the size $N$ of the moment set.
 
In turbulent flows, the choice of the closure law is directly related to the occurrence of Particle Trajectory Crossings (PTC).
When the Stokes number of the particles with respect to the Kolmogorov time scale $\text{St}_K=\tau_p/\tau_K$ is smaller than one, the velocity distribution is monokinetic i.e. the particulate phase has only one velocity per position in physical space. Several strategies exist in the literature for this range of Stokes number:
\begin{itemize}
\item the dusty gas approach \cite{carrier1958,marble1970}: the particles have the same velocity as the gas phase;
\item the Equilibrium Eulerian approach \cite{balachandar2010}: particle velocity is chosen has an expansion around the gas velocity;
\item the Monokinetic approach \cite{laurent2001}: the particle velocity is solved using a conservation equation for the momentum of the particulate phase. This approach had also been envisioned in a volume-averaged sense in \cite{druzhinin98}.
\end{itemize}

When $\text{St}_K>1$, the velocity distribution becomes multivalued because of PTC, referred as Random Uncorrelated Motion (RUM) in \cite{simonin2002}, and the closure laws for higher order moments has to be chosen to reproduce the NDF as accurately as possible. Moreover the choice of the NDF is not uniquely determined, and, as stated by the Hamburger Moment problem, for a finite set of moments, an infinite number of distributions is possible. So, two choices have to be made: (1) the number of moments to solve and (2) the assumptions on the unclosed moments or the shape of the NDF. In fact, this second choice splits the literature into two type of methods:
\begin{itemize}
\item Algebraic-Closure-Based moment methods (ACBMM)\cite{masi2014ijmf,masi2014ftc}\footnote{These groups were also referred as differential and algebraic model respectively in \cite{alipchenkov2007}}: The closure is devised based on the use of a limited information on the moments, like the total energy of the particulate phase, and  a series of assumptions on the moments, for example equilibrium assumption on the deviatoric part of the stress tensor, see \cite{kaufmann2008};
\item Kinetic-Based moment methods (KBMM) \cite{laurent12brief}: the closure is made at the kinetic level, by choosing a presumed equilibrium shape for the NDF.
\end{itemize}
The former are interesting, because it provides closed systems using a limited number of moments (density, momentum, and central energy). However, as the link between the moment and kinetic levels is not straightforward, it faces reliability issues i.e. the presumed moments do not correspond to a positive NDF. The KBMM overcome this issue by relying on a positive presumed shape for the NDF, thus realizability is intrinsically preserved. Moreover, as KBMM generates well-defined systems of equations of hyperbolic or weakly-hyperbolic type with an entropic structure \cite{Lev96,Lev98,kah11_cms,vie2013cicp}, the design of realizable numerical methods is easier than for ACBMM, for which the mathematical structure is hard to determine. Thus in this work we will focus on KBMM, and our choice will be driven by the needs of the LES modelling.


 \subsection{The Filtered Population Balance Equation}

In real turbulent applications, the range of scale encountered in the carrier phase is so wide that solving the whole spectrum is totally unreachable with the available computational resources. To avoid this problem, LES methods filter equations in physical or  frequency space. The filtering operation for a quantity $\phi$ is:
\begin{equation}
\overline{\phi}=\iiint G(\xx-\xx_1,\Delta)\phi(\xx_1,t)\text{d}\xx_1,
\end{equation}
where $\overline{\cdot}$ denotes a filtered quantity. The filtering of the carrier phase will decompose the velocity seen by the particle $\uu_g=\overline{\uu}_g+\uu_g'$  where  $\overline{\uu}_g$ is the filtered carrier phase velocity and $\uu_g'$ is the carrier phase velocity fluctuation.  
For the disperse phase, two methods have been used in the literature:
\begin{itemize}
\item Filtering the moment equation \ref{momenteq}, see \cite{shotorban2007,moreau2010};
\item Filtering the PBE and get moments of the filtered NDF, see \cite{pandya2002,zaichik2009}.
\end{itemize}
Theoretically the two methods lead to the same equations, as velocity and physical spaces are independent.
However, the latter strategy seems more interesting because:
\begin{itemize}
\item the PBE is linear whereas the moment equations are non-linear;
\item the filtering at the moment level loses track of the link with the underlying NDF.
\end{itemize}
Consequently, filtering at the kinetic level is our choice for the present work.

Filtering Eq.~(\ref{WBE}), we get:

\begin{equation}
\partial_t \overline{f}+\vv\cdot\partial_\xx \overline{f}+\partial_\vv\left( \dfrac{\overline{\uu}_g-\vv}{\taup} \overline{f}  \right)=-\dfrac{1}{\taup}\partial_\vv\cdot(\uu_g'f)^r,
\label{WBEugfilt}
\end{equation}
where $\overline{f}$ is the filtered NDF and $(\uu_gf)^r=\overline{\uu_gf}-\overline{\uu}_g\overline{f}$ is the subgrid correlation between the carrier phase velocity and the NDF of the particles.

To model this subgrid correlation, two contributions have been found in the literature that lead to the same closure:
\begin{itemize}
\item In \cite{pandya2002}, the authors derive the kinetic equation in direct analogy with the work of Reeks \cite{reeks91,reeks1992} who use the Lagrangian History Direct Interaction (LHDI);
\item In \cite{zaichik2009}, the authors follow the same strategy but based on the work of Zaichik et al. \cite{zaichik1999}, which consider the impact of the gas phase on the disperse phase to be a Gaussian process, and thus make use of the Furutsu-Donsker-Novikov Formula \cite{novikov65}. 
\end{itemize}

Following \cite{zaichik2009} , the closure for the subgrid correlation in the kinetic equation is:
\begin{equation}
-\dfrac{1}{\taup}\partial_\vv \cdot (\uu_g'f)^r=\partial_\vv\cdot \left( \boldsymbol \mu \partial_\vv \overline{f} + \boldsymbol \lambda \partial_\xx \overline{f} \right),
\label{zaichikclosure3D}
\end{equation}
where:
\begin{align}
&\lambda_{ij}=g_u^r\tau_{g,ij}, \quad \mu_{ij}=\tau_{g,ik}\left(\dfrac{f_u^r}{\tau_p}+l_u^r\partial_{x_k}\overline{u}_j \right),
\end{align}
where $f_u^r$, $g_u^r$ and $l_u^r$ are response coefficients (see \cite{zaichik2009} for details) whose long-time values are:
\begin{align}
f_u^r=\dfrac{1}{1+\St}, \quad g_u^r=\dfrac{1}{\St(1+\St)},\quad l_u^r =\dfrac{1}{\St(1+\St)^2}.
\end{align}
Here, the Stokes number is $\St=\tau_p/\tau_L^r$ where $\tau_L^r$ is the Integral time scale of the subgrid eddies.

The resulting kinetic equation in non-dimensional form is:
\begin{equation}
\partial_{t^*} \overline{f}+\vv^*\cdot\partial_{\xx^*} \overline{f}+\partial_{\vv^*}\cdot\left( \dfrac{\overline{\uu}_g^*-\vv^*}{\text{St}} \overline{f} \right)=\partial_{\vv^*}\cdot\left( \boldsymbol{\mu}^* \cdot\partial_{\vv^*} \overline{f} + \boldsymbol{\lambda}^* \cdot\partial_{\xx^*}\overline{f} \right),
\label{WBEfinal}
\end{equation}
where:
\begin{equation}
\mu_{ij}^*=\tau_{g,ik}^*\left(\dfrac{f_u^r}{\St}\delta_{kj}+ l_u^r \dfrac{\partial \overline{u}_{g,j}^*}{\partial x_k^*} \right), \quad \lambda_{ij}^*=g_u^r\tau_{g,ij}^*.
\end{equation}
In the following, star exponents are dropped for the sake of clarity.

\subsection{Asymptotic limit of the filtered PBE for small Stokes number}
To determine the equilibrium state and the asymptotic limit of the filtered kinetic equation Eq.~(\ref{WBEfinal}), The Chapman-Enskog expansion is used \cite{chapman1939}. Contrary to the one proposed in \cite{alipchenkov2007}, it is performed in the classical way, assuming a decomposition of the solution into power of $\St$: 
\begin{equation}
\tf=\tf^0+\St \tf^0\phi^1+\mathcal{O}\left(\St^2\right).
\label{CEdecomp}
\end{equation}
Rewriting Eq.~(\ref{WBEfinal})
\begin{equation}
\partial_{t} \overline{f}+\vv.\partial_{\xx} \overline{f}=\dfrac{1}{\St}\mathcal{J}(\tf),
\label{WBEreduce}
\end{equation}
where $\mathcal{J}(f)$ is the relaxation operator:
\begin{equation}
\mathcal{J}(\tf)=\partial_{\vv}\cdot\left[ (\vv-\overline{\uu}_g) \overline{f} + \St\boldsymbol{\mu}\partial_{\vv} \overline{f} + \St\boldsymbol{\lambda}\partial_{\xx}\overline{f}  \right].
\end{equation}
Inserting Eq.~(\ref{CEdecomp}) into Eq.~(\ref{WBEreduce}), and grouping terms by powers of $\St$:
\begin{equation}
\partial_{t} \tf^0+\vv\cdot \partial_{\xx} \tf^0+\St\left(\partial_{t} \tf^0\phi^1+\vv\cdot\partial_{\xx} \tf^0\phi^1 \right)=\dfrac{1}{\St}\mathcal{J}(\tf^0)+\mathcal{J}(\tf^0\phi^1)+\mathcal{O}(\St).
\label{CEequation}
\end{equation}

For small Stokes number, the zeroth order of  Eq.~(\ref{CEequation}) is:
\begin{equation}
\mathcal{J}(\tf^0)=\St\boldsymbol{\lambda} \partial_\xx\tf^0 +\St\boldsymbol{\mu}\partial_\vv\tf^0 - \left(\overline{\uu}_g-\vv\right) \overline{f}^0 =0.
\label{WBEasymptotic2}
\end{equation}
The solution of such an equation has the form:
\begin{align}
&\tf^0(x,u)=\Gamma(\xx,\vv)\dfrac{\rho}{(2\pi)^{3/2}\sqrt{|\boldsymbol\tau_g|}}\exp\left(-(\vv-\overline{\uu}_g)\boldsymbol \tau_g^{-1}(\vv-\overline{\uu}_g)^T\right) ,
\end{align}
where $\rho=\int \vv f \text{d}\vv$ and $\Gamma(\xx,\vv)$ is an unknown function that is equal to 1 in the case of an homogenous NDF in space.
The zeroth order distribution leads to the following equation on density:
\begin{align}
\partial_t \rho+\partial_\xx\cdot\rho \uu^0=0,
\end{align} 
where $\uu^0=\int \vv f^0\text{d}\vv$. To find $u^0$, the zeroth order moment of Eq.~(\ref{WBEasymptotic2}) is taken and leads to:
\begin{align}
\uu^0=\overline{\uu}_g-\dfrac{\boldsymbol\tau_g}{\rho}\partial_\xx\rho,
\end{align}
so that the asymptotic limit at the kinetic level leads to the following asymptotic limit at the moment level:
\begin{align}
\partial_t\rho+\partial_\xx\cdot\rho \overline{\uu}_g=\partial_\xx\cdot\left( \boldsymbol\tau_g\partial_\xx \rho \right).
\end{align}
One important thing to notice at this point is that this asymptotic limit is not a consequence of the modelling approach we choose, this is the limit that every LES model for the particulate phase has to recover in the low Stokes number limit. It is also worth noticing that the resulting moment system does not depend on any closure for the disperse phase, but only on the closure for the gas phase. Thus, any moment method that treats the number density should have the same asymptotic limit.
To achieve the description of the low Stokes number limit, taken the first order moment of Eq.~(\ref{WBEasymptotic2}), the asymptotic limit of the internal energy $\boldsymbol\sigma^0=\frac{1}{2}\int \vv\otimes\vv f^0\text{d}\vv-\frac{1}{2}\uu\otimes\uu$ is recovered:
\begin{align}
\boldsymbol\sigma^0=\boldsymbol\tau_{g}\left(1-\partial_\xx\otimes\overline{\uu}_g+\partial_\xx\otimes\left(\dfrac{\boldsymbol\tau_g}{\rho}\partial_\xx\rho \right)\right),
\end{align}
where $\otimes$ is the symmetric tensor outer product.

\subsection{Moment equations for LES}
After obtaining the kinetic equation, one has to go up to the moment level, and thus to choose the number of moments to solve and the closure law. Our choice is driven by the asymptotic limit of the PBE. Actually, the equilibrium distribution is a perturbation of a Gaussian distribution. Therefore, the minimal model to reproduce this asymptotic limit has to be a Gaussian distribution. Following \cite{vie2013cicp}, where the Anisotropic Gaussian distribution is used to close the moment system in a DNS context, we will use moments up to second order, that is 10 moments in 3D. The resulting system of equation is the following \footnote{In \cite{pandya2002}, the authors derive equations for zero-to-second order moments and  close the fluxes by using a zero third order central moments assumption, but they do not rely on the link with the kinetic level, see \cite{vie2013cicp} for details.}:
\begin{align}  \label{eq:AnisoGauss}
  \partial _t \rho+\partial _ { \xx} \cdot (\rho \uu) &=0 \\ 
  \partial _t (\rho  \uu)+\partial _ { \xx} \cdot (\rho\uu\uu^T+\mathbf{P}) &=\rho \dfrac {\uu_g-\uu}{\St} \\
 \partial _t (\rho \mathbf E)+\partial _ { \xx} \cdot ( (\rho \mathbf E + \mathbf P)\otimes \uu) &=\rho\dfrac { \overline{\uu}_g \otimes \uu +\St\boldsymbol{\mu}-2\mathbf E}{\St} ,
\end{align}
where  $E_{ij}=\frac{1}{2}\tilde{u}_i\tilde{u}_j+\frac{\sigma_{ij}}{2}$ is the total energy tensor, and $P_{ij}$ the pressure tensor:
\begin{align}
\rho\sigma_{ij}&=\iiint \left[(v_i-u_i)(v_j-u_j)\right]\overline{f}(t,\xx,\vv)\text{d}\vv \\
P_{ij}&=\overline{\rho}\left(\sigma_{ij}+\lambda_{ij}\right).
\end{align}

\subsection{Numerical issues associated to the low Stokes number limit}
In the following, the study will be fully performed in 1D for the sake of simplicity, even if all the developments can be straightforward envisaged in 3D. Moreover, we consider an homogeneous gaseous flow field, to focus on the main difficulties of the moment method. 

The resulting system of equation is:

\begin{equation} 
\left\{
\begin{array}{l}
\partial_t \rho+\partial_x \rho u = 0,\\
\partial_t \rho u+\partial_x\left( \rho u^2+P \right)=-\rho\dfrac{u-\overline{u}_g}{\text{St}},\\
\partial_t \rho E+\partial_x \left( \rho uE+uP\right)=-\rho u \dfrac{u-\overline{u}_g}{\text{St}} -\rho  \dfrac{2\epsilon-\text{St}\mu}{\text{St}}
\end{array}
\right.
\label{euler}
\end{equation}
where $E=u^2/2+\epsilon$ is the total energy and $P=\rho(2\epsilon+\lambda)$ is the pressure. 
Observe from now on that the last equation can be formulated in terms of internal energy and for smooth solutions as follows
$$
\partial_t \epsilon + u \partial_x \epsilon + \frac{P}{\rho} \partial_x u = 
- \dfrac{2\epsilon-\text{St}\mu}{\text{St}}.
$$
Recall that $\lambda$ 
and $\mu$ are given by 
\begin{equation} \label{lambda}
\mu=\dfrac{\tau_g}{\St(1+\St)},\quad \lambda=\dfrac{\tau_g}{\text{St}\left(1+\text{St} \right)}.
\end{equation}
Therefore, it is clear that the subgrid scale effects represented by $\tau_g \neq 0$ appear in two terms: a relaxation term on the internal energy, which tends to relax the internal energy of the disperse phase towards the subgrid scale internal energy of the gas phase (see paragraph below on the asymptotic analysis), and a pressure-like term via $\lambda$, acting like a flux and which propagates the time and space correlation of the gas phase over the disperse phase. 

{\it Characteristic analysis.} Easy manipulations show that 
(\ref{euler}) 
equivalently writes for smooth solutions 
\begin{equation} \label{modele_eulerien_primitive}
\left\{
\begin{array}{l}
\displaystyle \partial_t \rho+ u \partial_x \rho + \rho \partial_x u = 0\\
\displaystyle \partial_t u + \frac{(2\epsilon+\lambda)}{\rho} \partial_x \rho + u \partial_x u 
+ 2 \partial_x \epsilon = -\dfrac{u-\overline{u}_g}{\text{St}}\\
\displaystyle \partial_t \epsilon + \frac{P}{\rho} \partial_x u + u \partial_x \epsilon = 
- \dfrac{2\epsilon-\text{St}\mu}{\text{St}}\\
\end{array}
\right.
.
\end{equation}
Focusing on the convective part only (without source terms), the characteristic velocities 
are then given by
$$
\lambda_1 = u-c \,\, < \,\, \lambda_2=u \,\, < \,\, \lambda_3 = u + c,
$$ 
where the speed of sound $c$ is given by $c=\sqrt{6\epsilon+3\lambda}$. Observe that the square of the sound speed 
is made of two contributions, namely the classical one $6\epsilon=\gamma(\gamma-1)e$ associated with the perfect gas  
equation of state with adiabatic coefficient $\gamma=3$, and a new contribution involving  the subgrid stress of the gas 
phase $\tau_g$ via $\lambda$. Importantly, the presence of this new contribution ($\tau_g \neq 0$) 
makes the sound speed $c$ tend to infinity when the Stokes number $\text{St}$ goes to zero. From a physical point of view, 
it simply means that the smaller the Stokes number $\text{St}$ is, the faster the time and space correlations of the gas phase propagate over the dispersed phase. From a numerical point of view, this property is expected to give a strong constraint 
on the time step for a fully explicit in time Godunov-type method according to the well-known CFL stability condition.

{\it Asymptotic analysis.}
Here we consider the asymptotic behavior of the model when the Stokes number $\text{St}$ tends to $0$. 
%
Multiplying the last two equations of (\ref{modele_eulerien_primitive}) by $\text{St}\left(1+\text{St}\right)$ 
and letting $\text{St}$ go to $0$ gives 
\begin{align} \label{newnew}
 u&=\overline{u}_g-\dfrac{ \tau_g}{\rho}\partial_x \rho, \\
2 \epsilon&=\tau_g \Big( \partial_x\left(\dfrac{\tau_g}{\rho}\partial_x\rho\right) +1 \Big).
\end{align}
Combining these two relations with the first equation of (\ref{modele_eulerien_primitive}) (written in conservation 
form) then gives the following asymptotic limit 
\begin{align}
\partial_t \rho+\partial_x\rho \overline{u}_g=\partial_x \left(\tau_g\partial_x\rho\right),&   \label{eulerAsymptotic} \\
2 \epsilon=\tau_g \Big( \partial_x\left(\dfrac{\tau_g}{\rho}\partial_x\rho\right) +1 \Big). & \label{ndedlla}
\end{align}
The number density then evolves according to an advection-diffusion equation, which is coherent with the asymptotic limit at the kinetic level.
If one want to solve small Stokes number flows only, it is easy to directly solve Eq.~(\ref{eulerAsymptotic}). If one want to be able to solve a large range of Stokes numbers with the same method, Eq.~(\ref{euler}) is needed. However, in the low Stokes number limit, Eq.~(\ref{euler}) encounter two strong and limiting effects:
\begin{itemize}
\item the sound speed $c$ and the source terms tend to infinity: it leads to strong constraints for explicit numerical schemes, as well as high numerical diffusion for material waves $u$.
\item the numerical method has to recover the asymptotic limit at low Stokes number, which is not possible with "classical" schemes.
\end{itemize}
In the following, a numerical scheme is proposed to handle 
these two features of the chosen methodology.

\section{Numerical scheme}
The aim of this section 
is to describe a numerical scheme which (1) recovers the advection-diffusion asymptotic limit for the number density, referred as the asymptotic-preserving property  \cite{jin2,jinlevermore}
and (2) can get rid of the small time step imposed by the source term and the acoustic waves in a standard 
time-explicit Godunov-type method.
We rely on the approach proposed by Chalons, Girardin and Kokh in \cite{chalons2012ap} and based on three  ingredients:
\begin{itemize}
\item a Lagrange-Projection decomposition \cite{godlewski1996}, 
\cite{coqueletal2010};
\item a relaxation strategy first introduced by Suliciu \cite{suliciu1998} and Jin and Xin \cite{jinxin}, and then further studied by many authors for instance in \cite{coqueletal2001,chalonscoquel,chalons2008,bouchut, coqueletal2010};
\item a USI (Upwinding Sources at Interfaces) approach first 
initiated by Cargo and Le Roux \cite{cargoleroux}, Greenberg and Le Roux \cite{GreenbergLeRoux}, Gosse 
and Le Roux \cite{GosseLeRoux} 
(see also for related works in a wide range of applications \cite{klar}, \cite{jin3}, \cite{JinPareschiToscani}, \cite{gossetoscani3}, \cite{buetcordier}, \cite{berthonturpault}, \cite{CarilloGoudonLafitte}, \cite{DegondDeluzetVignalSangam}, \cite{FilbetJin}, \cite{DespresBuetFranck}, \cite{bouchutetal}, \cite{berthonleflochturpault,abgrallberthon}, \cite{HaackJinLiu} and the references therein). \\
\end{itemize}


The next three subsections give only a brief description of these three main ingredients as details 
can be found for instance in \cite{chalons2012ap}. The objective is to provide the reader with the most important 
update formulas and to focus on the specific treatment of the singular fluxes as well as of the source term associated with the internal energy, which are not included in the upwinding at the interfaces, receives a treatment through a splitting like approach and still satisfied the objective of being asymptotic preserving on energy. 

\subsection{Lagrange-Projection and source terms decomposition}
The objective of the Lagrange-projection is to decompose the full system into two sub-systems using a 
chain rule argument on the space derivatives. The first one only involves the transport wave, and the second 
one involves the acoustic waves (and the source terms). This kind of approach can be seen as an operator splitting 
strategy, and the main interest is to enable the use of different solvers for each subsystem. 
In this work, we shall consider explicit-explicit solvers but also implicit-explicit solvers 
in order to get rid of the strong CFL stability restriction imposed by the sound speed. By implicit-explicit, 
we mean here implicit on the acoustic waves and source terms and explicit on the transport part following the same approach as in \cite{chalons2012ap} (see \cite{coqueletal2010}). Note that the source 
term associated with the internal energy in the last equation of (\ref{euler}) being not considered in 
\cite{chalons2012ap}, we shall treat it 
separetely here using again a splitting strategy. Therefore, we shall end with three sub-systems to 
be treated numerically. \\
Using the property 
$\partial_x {\rho u X}= u \partial_x {\rho X} + \rho X \partial_x u$ for $X=\{1, u, E\}$ in the full 
system (\ref{euler}), we get the following transport system:
\begin{equation}
\left\{
\begin{array}{l}
\partial_t \rho+u\partial_x \rho =0 \\
\partial_t \rho u+u\partial_x\rho u=0 \\
\partial_t \rho E+u\partial_x \rho E=0,
\end{array}
\right.
\label{material}
\end{equation}
which corresponds to the advection of each field with velocity $u$. The acoustic and source term part 
is then given by
\begin{equation}
\left\{
\begin{array}{l}
\partial_t \rho+\rho\partial_x u =0 \\
\partial_t \rho u+\rho u\partial_x u+\partial_xP=-\rho\dfrac{u-\overline{u}_g}{\St}\\
\partial_t \rho E+\rho E\partial_x u +\partial_xPu=
-\rho u\dfrac{u-\overline{u}_g}{\St}-\rho\dfrac{2\epsilon-\St \mu}{\St}. 
\end{array}
\right.
\label{acoustic}
\end{equation}
Introducing $\tau=1/\rho$ and $\tau\partial_x=\partial_m$ this sytem equivalenty writes
\begin{equation}
\left\{
\begin{array}{l}
\partial_t \tau-\partial_m u =0 \\
\partial_t u+\partial_mP=-\dfrac{u-\overline{u}_g}{\St}\\
\partial_t E+\partial_mPu=-u\dfrac{u-\overline{u}_g}{\St}-\dfrac{2\epsilon-\St \mu}{\St}.
\end{array}
\right.
\label{acoussimp}
\end{equation}
Now splitting the drag force and energy relaxation terms leads to 
\begin{equation}
\left\{
\begin{array}{l}
\partial_t \tau-\partial_m u =0 \\
\partial_t u+\partial_mP=-\dfrac{u-\overline{u}_g}{\St}\\
\partial_t E+\partial_mPu=-u\dfrac{u-\overline{u}_g}{\St},
\end{array}
\right.
\label{acoussimpss}
\end{equation}
and
\begin{equation}
\left\{
\begin{array}{l}
\partial_t \tau=0 \\
\partial_t u=0\\
\partial_t \epsilon =-\dfrac{2\epsilon-\St \mu}{\St}.
\end{array}
\right.
\label{acoussimpssss}
\end{equation}
It is now a matter of approximating the three systems (\ref{acoussimpss}), (\ref{material}) and 
(\ref{acoussimpssss}).

\subsection{Relaxation procedure}
As mentioned above, system (\ref{acoussimpss}) may receive an implicit treatment in order to reduce 
the CFL stability restriction. In order for this implicit procedure to be cheap we propose to 
use a relaxation strategy, the objective of which being to overcome the difficulties introduced by the 
pressure nonlinearities. More precisely, the main idea is to consider the pressure 
$P$ as a new unknown $\Pi$ of the system and to solve for an additional equation for this new variable. 
We propose the following relaxation system
\begin{equation}
\left\{
\begin{array}{l}
\partial_t \tau-\partial_m u =0 \\
\partial_t u+\partial_m\Pi=-\dfrac{u-\overline{u}_g}{\St}\\
\partial_t E+\partial_m\Pi u=-u\dfrac{u-\overline{u}_g}{\St}  \\
\partial_t \Pi + a^2\partial_m u=0.
\end{array}
\right.
\label{relax}
\end{equation}
According to the very classical relaxation method, see for instance \cite{suliciu1998}, \cite{jinxin}, 
\cite{chalons2008}, \cite{bouchut}, this system 
will be solved instead of (\ref{acoussimpss}) at each time step and with $\Pi = P$ 
at initial time (the new relaxation variable $\Pi$ is said to be at equilibrium). 
It is worth noticing that the convective part of (\ref{relax}) has three constant eigenvalues given by 
$-a$, $0$ and $a$ so that the associated characteristic fields are linearly degenerate. Note also that
$a$ has to be chosen sufficiently large according to the subcharacteristic condition $a>\max\left(\rho c\right)$ 
in order to avoid instablities in this relaxation procedure. We refer for instance the reader to \cite{chalonscoquel} 
and \cite{coquelseguin} for a rigorous proof and for more 
details. \\
To conclude this short section, let us introduce the new variables $\overleftarrow{w}$ and 
$\overrightarrow{w}$ defined by $\overleftarrow{w}=\Pi - au$ and $\overrightarrow{w}=\Pi + au$. These 
quantities are easily shown to be strong Riemann invariants associated with the characteristic speeds $\pm a$ 
and allow to equivalently write (\ref{relax}) as follows~: 
\begin{equation}
\left\{
\begin{array}{rcl}
\partial_t \tau-\partial_m u &=&0\\
\partial_t \overrightarrow{w} + a\partial_m\overrightarrow{w}&=&-\dfrac{u-\overline{u}_g}{\St}\\
\partial_t \overleftarrow{w} -a\partial_m\overleftarrow{w}&=&\dfrac{u-\overline{u}_g}{\St}\\
\partial_t E+\partial_m\Pi u&=&-u\dfrac{u-\overline{u}_g}{\St}.  
\end{array}
\right.
\label{relaxsimp}
\end{equation}

\subsection{Notations and overview of the numerical scheme}

In this section, we introduce some notations and present the proposed numerical method 
for approximating the solutions of (\ref{euler}). This method can be understood as an
operator-splitting strategy consisting in solving the three systems (\ref{acoussimpss}), 
(\ref{acoussimpssss}) and (\ref{material}) in this order. Recall that (\ref{acoussimpss}) 
will be solved using the relaxation system (\ref{relaxsimp}). \\
Let be given a constant time step $\Delta t$ and a constant space step $\Delta x$. Introducing
$x_{j+1/2}=(j+1/2) \Delta x$ for $j \in \mathbb{Z}$ and $t^n=n \Delta t$ for $n \in \mathbb{N}$, 
the approximate solution of (\ref{euler}), say ${\bf u}_{\lambda}(x,t)$, with ${\bf u}_0$ as initial
data is classically sought as a piecewise constant function on each slab $C^n_j=[x_{j-1/2};x_{j+1/2})
\times [t^n;t^{n+1})$ for $(n,j) \in \mathbb{N} \times \mathbb{Z}$, and we set
$$
{\bf u}_{\lambda}(x,t) = {\bf u}^n_j \quad \mbox{for all} \quad (x,t) \in C^n_j, \quad j \in \mathbb{Z},
\quad n \in \mathbb{N}.
$$
Here $\lambda$ stands for the ratio $\lambda= \Delta t / \Delta x$, and for completeness we set
$$
{\bf u}^0_j = \frac{1}{\Delta x} \int_{x_{j-1/2}}^{x_{j+1/2}} {\bf u}_0(x) dx.
$$
Assuming that the approximate solution ${\bf u}_{\lambda}(x,t^n)$ is known, we propose to
advance it to the next time level $t^{n+1}$ according to a classical splitting
strategy between (\ref{acoussimpss}), 
(\ref{acoussimpssss}) and (\ref{material}). 
%
\ \\
\ \\
{\bf First step} ($t^{n} \to t^{n+1=}$) This step consists in solving the Lagrangian system (\ref{acoussimpss}), 
that is to say
\begin{equation} \label{sifrac2}
\partial_t {\bf v} + \partial_m {\bf g(v)} = {\bf S(v)},
\end{equation}
with ${\bf v}=(\tau,u,E)^t$, ${\bf g(v)}=(-u,P,Pu)^t$ and 
${\bf S(v)}=(0,-\dfrac{u-\overline{u}_g}{\St},-u\dfrac{u-\overline{u}_g}{\St})^t$. Note that the change of variables 
${\bf u} \to {\bf v}$ is one-to-one. We consider 
${\bf v}_{\lambda}(x,t^{n})$ (easily calculated from ${\bf u}_{\lambda}(x,t^{n})$) as initial 
data. 
To define $\bf{u}^{n+1=}_j$, we propose to use a Godunov-type method based on the definition 
of an approximate Riemann solver for (\ref{relaxsimp}) as detailed in \cite{chalons2012ap}. The proposed approximate Riemann solver here coincides exactly with the one 
proposed in \cite{chalons2012ap} up to very minor modifications~: the gravity is not taken into account in the present paper  and 
the friction terms make the velocity tend to $\overline{u}_g$ here instead of $0$.
This approximate Riemann solver includes the source terms in its definition in the sense that the definition 
of the intermediade states actually depends on the source term ${\bf S(v)}$ and in particular on 
the Stokes number $\text{St}$. Therefore, the source term ${\bf S(v)}$ is taken into account at each interface of the mesh, 
in agreement with the celebrated USI (Upwinding Sources at Interfaces) approach introduced in 
\cite{cargoleroux}, \cite{GreenbergLeRoux}, \cite{GosseLeRoux}. Note however that the way the sources are upwinding here is really different 
and relies on the formalism of {\it consistency in the integral sense} introduced by Gallice in \cite{gallice2003b}, \cite{gallice2003} for systems of conservation 
laws with sources. This formalism is nothing but a generalization of the well-known Harten, Lax and van Leer 
formalism \cite{hlvl} for systems of conservation laws. We then refer the reader to \cite{chalons2012ap} for more details. \\
\ \\
{\it Time-explicit Godunov-type scheme.}
Averaging the approximate Riemann solutions defined at each interface leads to the following Godunov-type 
update formulas 
\begin{equation}\label{GIRpart2syst11}
\left\{ 
\begin{aligned}
\tau _{j}^{n+1=} 
&=
 \tau_{j}^{n} + \frac{\Delta t}{\Delta m_j}(u_{j+\frac{1}{2}}^{*} - u_{j-\frac{1}{2}}^{*})
,\\
\overrightarrow{w}_{j}^{n+1=} 
&=
 \overrightarrow{w}_{j}^{n} 
- a\frac{\Delta t}{\Delta m_j}(\overrightarrow{w}_{j}^{n}-\overrightarrow{w}_{j-1}^{n}) + 
\Delta t \, a \, \frac{\Delta m_{j-1/2}}{\Delta m_{j}} \frac{\overline{u}_g-u^{*}_{j-\frac{1}{2}}}{St}
,\\
\overleftarrow{w}_{j}^{n+1=} 
&=
 \overleftarrow{w}_{j}^{n} 
+ a\frac{\Delta t}{\Delta m_j}(\overleftarrow{w}_{j+1}^{n}-\overleftarrow{w}_{j}^{n}) -  
\Delta t \, a \, \frac{\Delta m_{j+1/2}}{\Delta m_{j}} \frac{\overline{u}_g - u^{*}_{j+\frac{1}{2}}}{St}
,\\
E_{j}^{n+1=} 
&=
 E_{j}^{n} 
- \frac{\Delta t}{\Delta m_j}((up)_{j+\frac{1}{2}}^{*}-(up)_{j-\frac{1}{2}}^{*}) 
 + 
\frac{\Delta t}{\Delta m_j} \overline{u}_g \frac{\Delta m_{j+1/2} u_{j+\frac{1}{2}}^{*} +
 \Delta m_{j-1/2} u_{j-\frac{1}{2}}^{*}}{2 St} \\
 &
- 
\frac{\Delta t}{\Delta m_j} 
\frac{\Delta m_{j+1/2} (u_{j+\frac{1}{2}}^{*})^{2}+
\Delta m_{j-1/2} (u_{j-\frac{1}{2}}^{*})^{2}}{2 St},
\end{aligned}\right.
\end{equation}
where we have set
$$
\Delta m_j =\rho^n_j \Delta x
,\quad
\Delta m_{j+1/2} = \frac{\Delta m_j + \Delta m_{j+1}}{2},
$$
and
\begin{equation} \label{ymnew}
u^{*}_{j+\frac{1}{2}} = \frac{St}{ 2 a St + \Delta m_{j+1/2}} 
( \overrightarrow{w}_{j}^{n} - \overleftarrow{w}_{j+1}^{n} + \overline{u}_g \frac{\Delta m_{j+1/2}}{St} ),
\quad 
p_{j+\frac{1}{2}}^{*} = \frac{\overrightarrow{w}_{j}^{n} + \overleftarrow{w}_{j+1}^{n}}{2}.
\end{equation}
Using the relation
$$
{u}_{j}^{n+1=} = \frac{\overrightarrow{w}_{j}^{n+1=} -\overleftarrow{w}_{j}^{n+1=} }{2a},
$$
we also have
\begin{equation}\label{GIRentropsysequ1}
\left\{ 
\begin{aligned}
\tau_{j}^{n+1=} 
=& \tau_{j}^{n} + \frac{\Delta t}{\Delta m_j}(u_{j+\frac{1}{2}}^{*} - u_{j-\frac{1}{2}}^{*}),\\
u_{j}^{n+1=} 
=& u_{j}^{n} 
-
\frac{\Delta t}{\Delta m_j}(p_{j+\frac{1}{2}}^{*} - p_{j-\frac{1}{2}}^{*}) + 
\frac{\Delta t}{\Delta m_j} 
\frac{
\Delta m_{j-1/2} (\overline{u}_g - u_{j-\frac{1}{2}}^*)
+ 
\Delta m_{j+1/2} (\overline{u}_g - u_{j+\frac{1}{2}}^*)
}{2 St}
,\\
E_{j}^{n+1=} =
&
E_{j}^{n} 
- \frac{\Delta t}{\Delta m_j}((up)_{j+\frac{1}{2}}^{*}-(up)_{j-\frac{1}{2}}^{*}) 
 + 
\frac{\Delta t}{\Delta m_j} \overline{u}_g \frac{\Delta m_{j+1/2} u_{j+\frac{1}{2}}^{*} +
 \Delta m_{j-1/2} u_{j-\frac{1}{2}}^{*}}{2 St} \\
 &
- 
\frac{\Delta t}{\Delta m_j} 
\frac{\Delta m_{j+1/2} (u_{j+\frac{1}{2}}^{*})^{2}+
\Delta m_{j-1/2} (u_{j-\frac{1}{2}}^{*})^{2}}{2 St}.\\
\end{aligned}\right.
\end{equation}
This scheme is shown to be stable under the following CFL condition
$$
\max_{j \in \mathbb{Z}} a \frac{\Delta t}{\Delta m_j} \leq \frac{1}{2},
$$
see again for instance \cite{chalons2012ap}. Using classical notations, one has thus defined 
${\bf v}_{\lambda}(x,t^{n+1=})$. This piecewise constant solution will be used in the second step as a natural 
initial condition. \\
\ \\
{\it Time-implicit Godunov-type scheme.}
In order to get rid of the above CFL restriction on the time step $\Delta t$, which becomes stronger and stronger 
as the Stokes number goes to zero, a time-implicit version can be simply defined as follows,
\begin{equation}\label{GIRpart2syst11imp}
\left\{ 
\begin{aligned}
\tau _{j}^{n+1=} 
&=
 \tau_{j}^{n} + \frac{\Delta t}{\Delta m_j}(u_{j+\frac{1}{2}}^{*} - u_{j-\frac{1}{2}}^{*})
,\\
\overrightarrow{w}_{j}^{n+1=} 
&=
 \overrightarrow{w}_{j}^{n} 
- a\frac{\Delta t}{\Delta m_j}(\overrightarrow{w}_{j}^{n+1=}-\overrightarrow{w}_{j-1}^{n+1=}) + 
\Delta t \, a \, \frac{\Delta m_{j-1/2}}{\Delta m_{j}} \frac{\overline{u}_g-u^{*}_{j-\frac{1}{2}}}{St}
,\\
\overleftarrow{w}_{j}^{n+1=} 
&=
 \overleftarrow{w}_{j}^{n} 
+ a\frac{\Delta t}{\Delta m_j}(\overleftarrow{w}_{j+1}^{n+1=}-\overleftarrow{w}_{j}^{n+1=}) -  
\Delta t \, a \, \frac{\Delta m_{j+1/2}}{\Delta m_{j}} \frac{\overline{u}_g - u^{*}_{j+\frac{1}{2}}}{St}
,\\
E_{j}^{n+1=} 
&=
 E_{j}^{n} 
- \frac{\Delta t}{\Delta m_j}((up)_{j+\frac{1}{2}}^{*}-(up)_{j-\frac{1}{2}}^{*}) 
 + 
\frac{\Delta t}{\Delta m_j} \overline{u}_g \frac{\Delta m_{j+1/2} u_{j+\frac{1}{2}}^{*} +
 \Delta m_{j-1/2} u_{j-\frac{1}{2}}^{*}}{2 St} \\
 &
- 
\frac{\Delta t}{\Delta m_j} 
\frac{\Delta m_{j+1/2} (u_{j+\frac{1}{2}}^{*})^{2}+
\Delta m_{j-1/2} (u_{j-\frac{1}{2}}^{*})^{2}}{2 St},
\end{aligned}\right.
\end{equation}
where $u^{*}_{j+\frac{1}{2}}$ and $p_{j+\frac{1}{2}}^{*}$ are now implicitly defined by
\begin{equation} \label{ymnewimp}
u^{*}_{j+\frac{1}{2}} = \frac{St}{ 2 a St + \Delta m_{j+1/2}} 
( \overrightarrow{w}_{j}^{n+1=} - \overleftarrow{w}_{j+1}^{n+1=} + \overline{u}_g \frac{\Delta m_{j+1/2}}{St} ),
\quad 
p_{j+\frac{1}{2}}^{*} = \frac{\overrightarrow{w}_{j}^{n+1=} + \overleftarrow{w}_{j+1}^{n+1=}}{2}.
\end{equation}
Let us notice that the second and third equations do not depend on $\tau$ and $E$ and can thus be solved 
independenlty. The corresponding system is easily seen to be linear with a pentadiagonal and strictly  
diagonally dominant matrix. Therefore, $\overleftarrow{w}_{j}^{n+1=}$ and $\overrightarrow{w}_{j}^{n+1=}$ 
are uniquely defined for any $\Delta t$. Then, $\tau_{j}^{n+1=}$ and $E_{j}^{n+1=}$ follow {\it explicitly} 
(which makes the overall strategy very low cost) 
thanks to the first and fourth equations of (\ref{GIRpart2syst11imp}). \\
\ \\
{\it About the Asymptotic-Preserving property.} 
Let us discuss in this short paragraph the asymptotic behavior of the proposed numerical scheme 
when the Stokes number $\text{St}$ goes to zero. We focus on the explicit-explicit version for simplicity. Let us 
first observe that  $\lambda St$ goes to $\tau_g$ when $\text{St}$ goes to zero by definition of $\lambda$, so that 
$St P$ goes to $\tau_g \rho$ by definition of $P$. We then easily get by definition of 
$\overleftarrow{w}$ and $\overrightarrow{w}$ the following limit
$$
\lim_{St \to 0} u^{*}_{j+\frac{1}{2}} = 
\frac{1}{ \Delta m_{j+1/2}} 
( \tau_g \rho_{j}^{n} - \tau_g \rho_{j+1}^{n} + \overline{u}_g {\Delta m_{j+1/2}} ),
$$
that is to say
$$
\lim_{St \to 0} u^{*}_{j+\frac{1}{2}} = \overline{u}_g - \frac{2 \tau_g}{\rho_{j}^n+\rho_{j+1}^n} 
\frac{\rho_{j+1}^{n}-\rho_{j}^{n}}{\Delta x}.
$$
This limit is clearly consistent with the expected mass flux 
$$
u=\overline{u}_g-\dfrac{ \tau_g}{\rho}\partial_x \rho,
$$
in (\ref{newnew}). \\ 
Regarding the asymptotic behavior of the internal energy (\ref{ndedlla}), it will be useful in the 
next steps of the method to estimate the asymptotic behavior when the Stokes number goes to zero 
of the following discrete time derivative 
\begin{equation} \label{tdae}
St \frac{\epsilon_j^{n+1=}-\epsilon_j^{n}}{\Delta t}.
\end{equation}
With this in mind, it is first easily shown from the second equation of (\ref{GIRentropsysequ1}) 
that $St u_j^{n+1=}$ goes to zero with the Stokes number (the calculations are left to the reader). 
Then letting $\text{St}$ go to zero in the third equation of (\ref{GIRentropsysequ1}) shows that (\ref{tdae}) 
asymptotically behaves like 
$$
St \Big( - \frac{\Delta t}{\Delta m_j}((up)_{j+\frac{1}{2}}^{*}-(up)_{j-\frac{1}{2}}^{*}) 
 + 
\frac{\Delta t}{\Delta m_j} \overline{u}_g \frac{\Delta m_{j+1/2} u_{j+\frac{1}{2}}^{*} +
 \Delta m_{j-1/2} u_{j-\frac{1}{2}}^{*}}{2 St} \\
$$
$$
- 
\frac{\Delta t}{\Delta m_j} 
\frac{\Delta m_{j+1/2} (u_{j+\frac{1}{2}}^{*})^{2}+
\Delta m_{j-1/2} (u_{j-\frac{1}{2}}^{*})^{2}}{2 St} \Big).
$$
Using the properties that $u_{j+\frac{1}{2}}^{*}$ and $St p_{j+\frac{1}{2}}^{*}$ are respectively 
consistent with $u=\overline{u}_g-\dfrac{ \tau_g}{\rho}\partial_x \rho$ and $\tau_g \rho$, together with the Leibniz relation
$$
(up)_{j+\frac{1}{2}}^{*}-(up)_{j-\frac{1}{2}}^{*} = 
\frac{1}{2}\big(  u_{j+\frac{1}{2}}^{*}+u_{j-\frac{1}{2}}^{*}  \big) 
\big( p_{j+\frac{1}{2}}^{*}-p_{j-\frac{1}{2}}^{*}   \big)
+
\frac{1}{2}\big(  p_{j+\frac{1}{2}}^{*}-p_{j-\frac{1}{2}}^{*}  \big) 
\big(  u_{j+\frac{1}{2}}^{*}-u_{j-\frac{1}{2}}^{*}  \big),
$$ 
we easily get that (\ref{tdae}) is asymptotically consistent with $- \tau_g \partial_x u$ 
when $\text{St}$ goes to zero. \\
\ \\
{\bf Second step} ($t^{n+1=} \to t^{n+1-}$) 
The second step consists in solving the transport step (\ref{material}) with 
${\bf v}_{\lambda}(x,t^{n+1=})$ as initial 
data, which can be easily and equivalently transformed 
in terms of the ${\bf u}=(\rho, \rho u, \rho E)^t$ variables.
Following \cite{godlewski1996} (see also \cite{chalons2012ap}), we consider a very classic upwind and time-explicit 
numerical scheme given by
\begin{align} \label{uflagstep}
{X}_j^{n+1-}={X}_j^{n+1=}+\dfrac{\Delta t}{\Delta x}\left[ u_{j-1/2}^{*,+} {X}_{j-1}^{n+1=} +\left[u_{j+1/2}^{*,-}-u_{j-1/2}^{*,+}\right] {X}_{j}^{n+1=}  - u_{j+1/2}^{*,-} {X}_{j+1}^{n+1=} \right] 
\end{align}
where $\alpha^{\pm}=(\alpha\pm|\alpha|)/2$ for any $\alpha$ and ${X}\in \{\rho,\rho u, \rho E\}$. This scheme is shown to be stable under the following CFL condition
\begin{equation} \label{ccflplalp}
\frac{\Delta t}{\Delta x} \big( u_{j-1/2}^{*,+} - u_{j+1/2}^{*,-} \big) \leq 1.
\end{equation}
\ \\
{\it About the Asymptotic-Preserving property.} Let us go on with the asymptotic-preserving property. The first equation
of (\ref{uflagstep}) writes 
$$
{\rho}_j^{n+1-}={\rho}_j^{n+1=}+\dfrac{\Delta t}{\Delta x}\left[ u_{j-1/2}^{*,+} {\rho}_{j-1}^{n+1=} +\left[u_{j+1/2}^{*,-}-u_{j-1/2}^{*,+}\right] {\rho}_{j}^{n+1=}  - u_{j+1/2}^{*,-} {\rho}_{j+1}^{n+1=} \right],
$$
that we can combine with the first equation of (\ref{GIRpart2syst11}), namely
$$
\tau _{j}^{n+1=} 
=
\tau_{j}^{n} \Big(1 + \frac{\Delta t}{\Delta x}(u_{j+\frac{1}{2}}^{*} - u_{j-\frac{1}{2}}^{*})\Big),
$$ 
to give the conservative update formula
\begin{equation} \label{cesr}
\rho _{j}^{n+1-} 
=
\rho_{j}^{n} - \frac{\Delta t}{\Delta x} \big( \{\rho u\}_{j+\frac{1}{2}} - \{\rho u\}_{j-\frac{1}{2}} \big),
\end{equation}
with numerical flux given by
$$
\{\rho u\}_{j+\frac{1}{2}} = \rho _{j}^{n+1=} u_{j+1/2}^{*,+} + \rho _{j+1}^{n+1=} u_{j+1/2}^{*,-}.
$$
We have just seen in the previous step that $u^{*}_{j+\frac{1}{2}}$ is consistent with  
$\overline{u}_g-\dfrac{ \tau_g}{\rho}\partial_x \rho$ in the limit $St \to 0$. As a consequence, 
the conservative formula (\ref{cesr}) is clearly consistent with 
the first equation (\ref{eulerAsymptotic}) of the expected asymptotic limit. \\Regarding the second equation 
(\ref{ndedlla}), we just note here that the following discrete time derivative
\begin{equation} \label{tdae2}
St \frac{\epsilon_j^{n+1-}-\epsilon_j^{n+1=}}{\Delta t}
\end{equation}
goes to zero when the Stokes number goes to zero since $\epsilon$ evolves according to 
the transport equation 
$$
\partial_t \epsilon + u \partial_x \epsilon = 0,
$$
in this step. Recall indeed that the interfacial velocity $u^{*}_{j+\frac{1}{2}}$ used to discretize this equation 
is such that $St u^{*}_{j+\frac{1}{2}}$ goes to zero with the Stokes number. \\
\ \\
{\bf Third step} ($t^{n+1-} \to t^{n+1}$) 
The last step consists in solving (\ref{acoussimpssss}) with 
${\bf v}_{\lambda}(x,t^{n+1-})$ as initial 
data.  Compared to the work of \cite{chalons2012ap}, this step is the main difference, as no internal energy relaxation exists in their work. Here we choose to decouple this relaxation step because the leading order for the asymptotic limit is on the mean velocity component, such as it is not mandatory to include this step in the HLLC solver. 
To define ${\bf v}_{\lambda}(x,t^{n+1})$, it is a natural idea to use 
a classical pointwise implicit evaluation of the source term. More precisely, 
it amounts to set
$$
\left\{
\begin{array}{l}
\tau^{n+1}_j=\tau^{n+1-}_j, \\
u^{n+1}_j=u^{n+1-}_j,\\
\epsilon^{n+1}_j =\epsilon^{n+1-}_j - \Delta t \dfrac{2\epsilon^{n+1}_j-\St \mu}{\St},
\end{array}
\right.
$$
or equivalently 
\begin{equation} \label{sifrac1}
\left\{
\begin{array}{l}
\tau^{n+1}_j=\tau^{n+1-}_j, \\
u^{n+1}_j=u^{n+1-}_j,\\
\epsilon^{n+1}_j =\dfrac{\St (\mu \Delta t + \epsilon^{n+1-}_j)}{\St+2 \Delta t}.
\end{array}
\right.
\end{equation}
{\it About the Asymptotic-Preserving property.} 
Let us first recall that the first equation (\ref{eulerAsymptotic})
has been proved at the end of the second step and remains valid since $\rho$ 
is not expected to be modified in this last step. Regarding now 
the second equation 
(\ref{ndedlla}), we first write (\ref{sifrac1}) under the following equivalent form
$$
St \frac{\epsilon^{n+1}_j -\epsilon^{n}_j}{\Delta t} + 2\epsilon^{n+1}_j =  
St \frac{\epsilon^{n+1=}_j -\epsilon^{n}_j}{\Delta t} + 
St \frac{\epsilon^{n+1-}_j -\epsilon^{n+1=}_j}{\Delta t} +\St \mu.
$$
Since 
$
\lim_{\St \to 0} \St \mu = \tau_g (1+\partial_x \overline{u}_g),
$
using the results obtained in the previous two steps then clearly gives that  
$2 \epsilon^{n+1}_j$ is asymptotically consistant with $-\tau_g \partial_x u + \tau_g + \tau_g \partial_x \overline{u}_g$, that 
is to say with (\ref{ndedlla}) by the relation $u=\overline{u}_g-\dfrac{ \tau_g}{\rho}\partial_x \rho$.

\subsection{Main properties}

We gather in this section the main properties of the proposed algorithm. We focus on the 
implicit-explicit version. \\

\begin{theorem}
\label{GIRpart3thrm1}
Under the CFL condition (\ref{ccflplalp}) and provided that $a$ is chosen sufficiently large, the implicit-explicit in time numerical scheme is well defined and 
satisfies the following stability properties:

{\it (i)} it is a conservative scheme for the density $\rho$. It is also a conservative scheme 
for $\rho u$ and $\rho E$ when the source terms are omitted,

{\it (ii)} the density $\rho^n_j$ is positive for all $j$ and $n>0$ provided that $\rho^0_j$ is positive for all $j$,

{\it (iii)} it is asymptotic preserving. \\
\end{theorem} 

{\it Sketch of the proof.} Property {\it (i)} has been proved in the course of the description of the second step 
for the density, see (\ref{cesr}). The proof is similar for $\rho u$ and $\rho E$. Property {\it (ii)} is 
obtained from standard manipulations \cite{godlewski1996}, \cite{coqueletal2010}. The asymptotic-preserving property has been 
proved in the previous subsections. \\
\ \\
{\it Remark.} The validity of an entropy inequality has been proved in \cite{chalons2012ap} for the explicit-explicit 
version of the scheme. Regarding the implicit-explicit version, we refer the reader to \cite{coqueletal2010}.

\section{Numerical results}
Here we propose a test case which highlights the effects of the subgrid scale of turbulence on the disperse phase. The domain is $[-1 ,1]$. The initial state for the disperse phase is a spatially Gaussian distribution $\rho(t=0,x)=1+\exp\left(-x^2/(2\sigma_0^2)\right)$, where $\sigma=0.01$, which mimics the dispersion 
of particles occurring in a turbulent field. At time $t=0$, the particles are at rest, i.e. $u(t=0,x)=0$ and $\epsilon(t=0,x)=0$. The gas phase is constant in time and space, and is decomposed into a mean velocity $\overline{u}_g=0$, which in fact generates no fluxes, and a subgrid scale energy $\tau_g=1/10$. 

For the explicit schemes, two different explicit constraints are imposed on the time step through the source terms and the CFL number:
\begin{equation}
\Delta t_{source}\le \dfrac{\St}{2}, \quad \Delta t_{CFL} \le \text{CFL}\dfrac{\Delta x}{c}=\text{CFL}\dfrac{\Delta x}{\sqrt{6\epsilon+3\lambda}}.
\end{equation} 
At time $t=0s$, the time step is then:
\begin{equation}
\Delta t=\min \left(\dfrac{\St}{2}, \dfrac{\text{CFL}}{\sqrt{3\tau_g}} \Delta x\sqrt{\St(1+\St)}        \right).
\end{equation}

%

\begin{figure}

\centering
\includegraphics[width=1.0\textwidth]{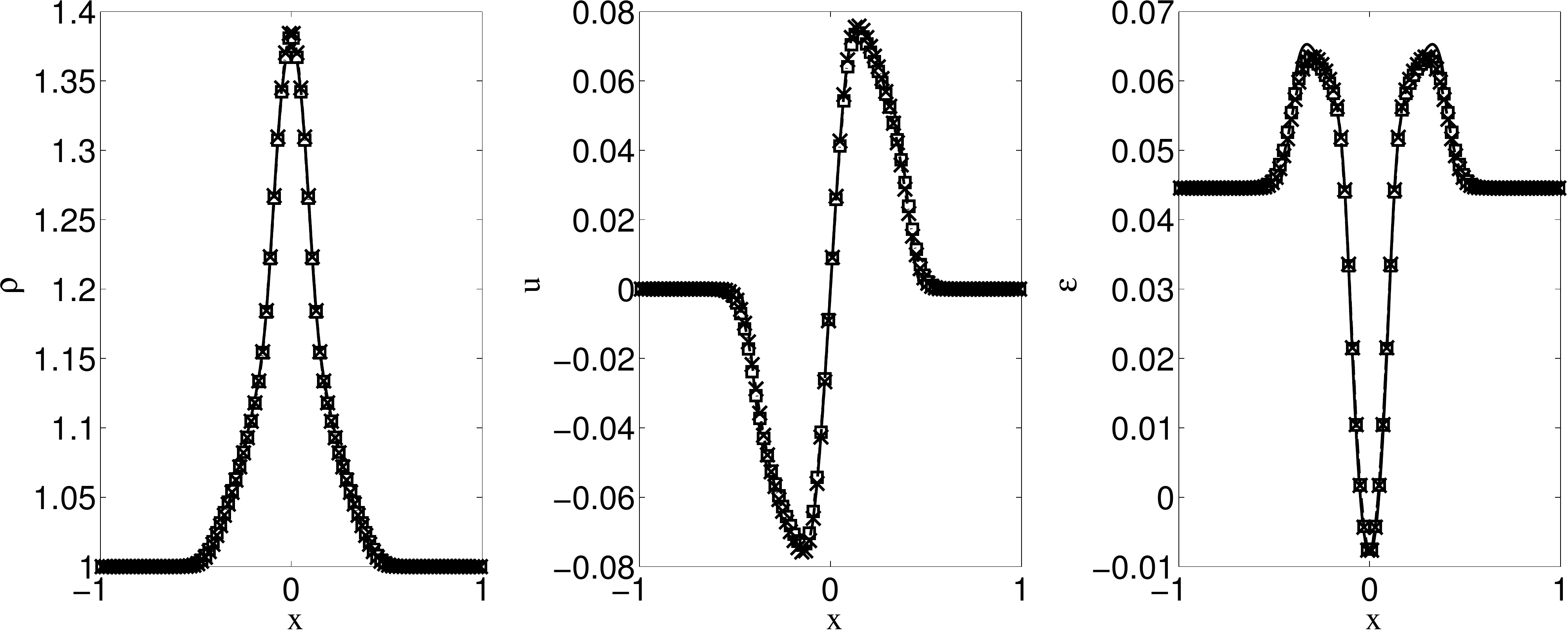}
\caption{Stokes number=0.1: Comparison between refined solution with 2000 cells (black full line), Explicit non-AP scheme (squares) and Explicit AP scheme (crosses) with 100 cells at time $t=0.2s$.}
\label{St0.1}
\end{figure}

\begin{figure}

\centering
\includegraphics[width=1.0\textwidth]{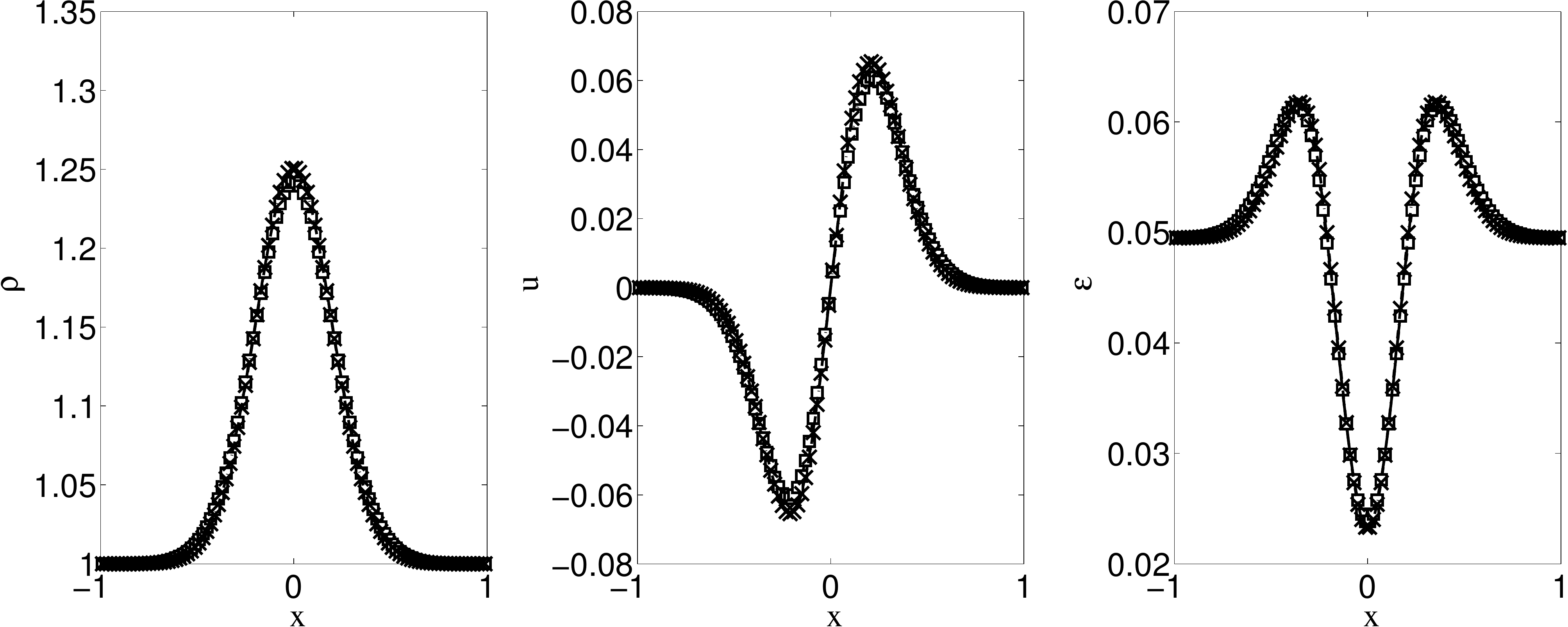}
\caption{Stokes number=0.01: Comparison between refined solution with 2000 cells (black full line), Explicit non-AP scheme (squares) and Explicit AP scheme (crosses) with 100 cells at time $t=0.2s$.}
\label{St0.01}
\end{figure}

In Figs.~\ref{St0.1}-\ref{St0.01}, results at time $t=0.2s$ for the density, the mean velocity and the internal energy are plotted against the position, for the Explicit non-AP and AP schemes, and for $\St=0.1$ and $\St=0.01$, using  100 cells.
The Stokes number is not sufficiently small to use the asymptotic solution as a reference. So a 2000-cell solution is use as a reference, for which a sufficient mesh convergence has been verified for the sake of the comparisons. For these Stokes number and this space discretization, the differences between non-AP and AP schemes are not obvious, as the contribution of the subgrid flux and source terms are not predominant.

\begin{figure}

\centering
\includegraphics[width=1.0\textwidth]{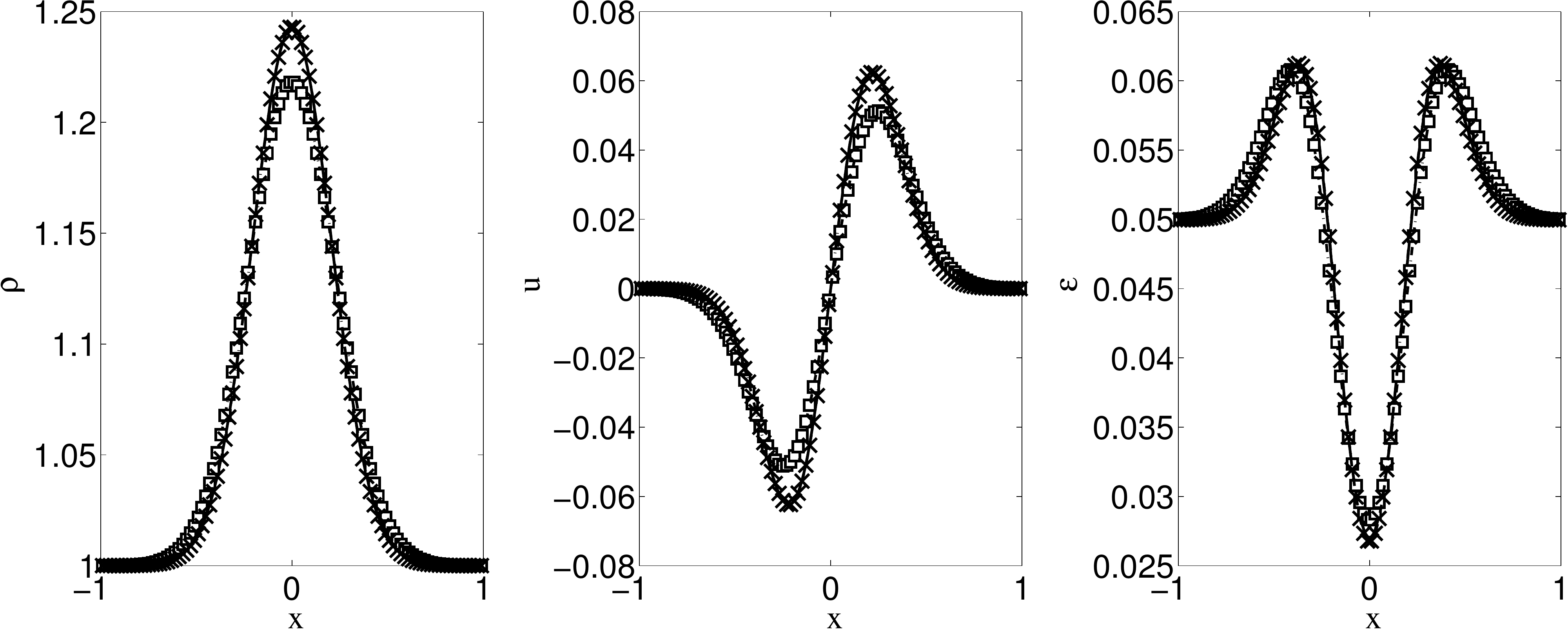}
\caption{Stokes number=0.001: Comparison between asymptotic analytic solution (black full line), Explicit non-AP scheme (squares) and Explicit AP scheme (crosses) with 100 cells at time $t=0.2s$.}
\label{St0.001}
\end{figure}

\begin{figure}

\centering
\includegraphics[width=1.0\textwidth]{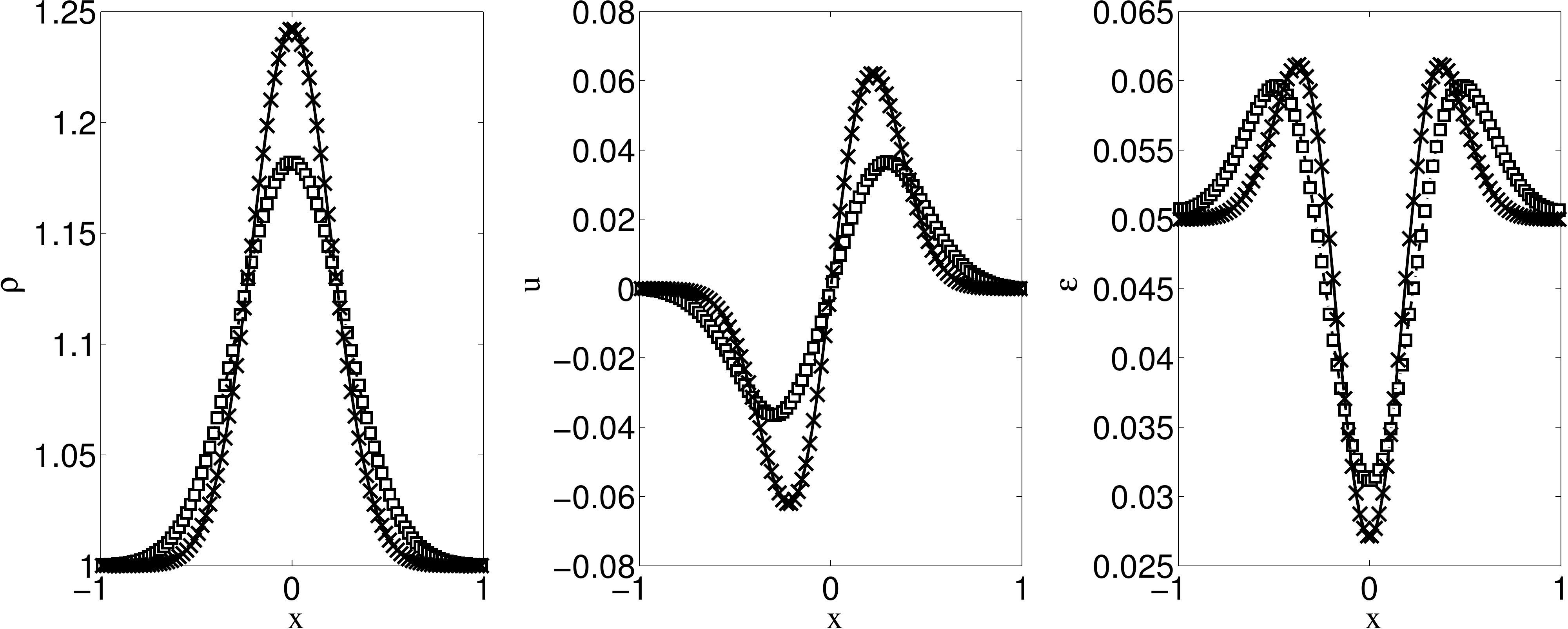}
\caption{Stokes number=0.0001: Comparison between asymptotic analytic solution (black full line), Explicit non-AP scheme (squares) and Explicit AP scheme (crosses) with 100 cells at time $t=0.2s$.}
\label{St0.0001}
\end{figure}

\begin{figure}

\centering
\includegraphics[width=0.5\textwidth]{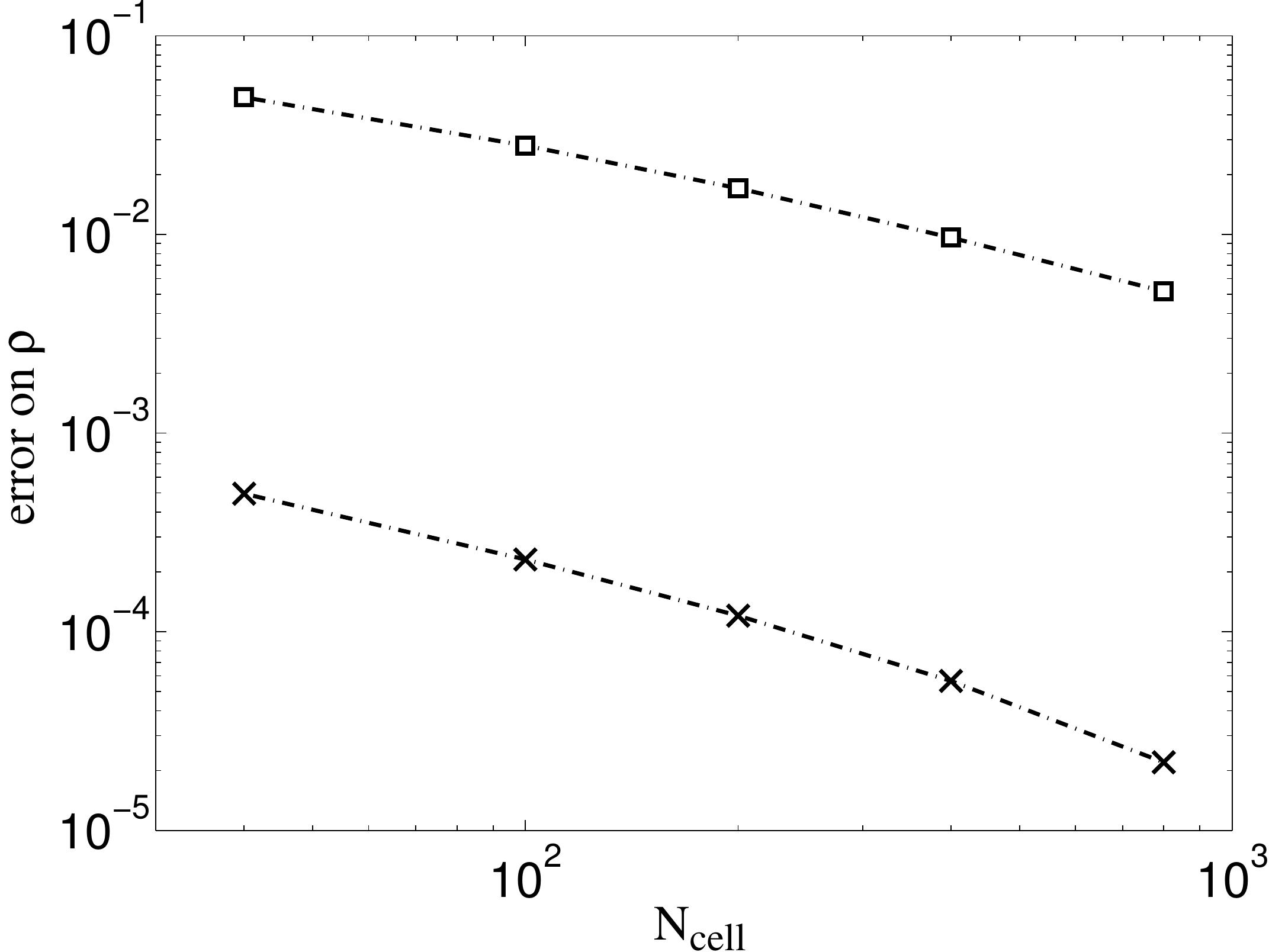}
\caption{Stokes number=0.0001: Error on the number density comparing to the asymptotic analytic solution using the Explicit non-AP scheme (squares) and the Explicit AP scheme (crosses).}
\label{St0.0001_ExplicitError}
\end{figure}

In Figs.~\ref{St0.001}-\ref{St0.0001}, results at time $t=0.2s$ for the density, the mean velocity and the internal energy are plotted against the position, for the Explicit non-AP and AP schemes, and for $\St=0.001$ and $\St=0.0001$, using  100 cells. Here the Stokes numbers are small enough for the time step to be CFL-driven. The Stokes number is now sufficiently small to use the asymptotic solution as a reference, the statement having been verified. Now the differences between non-AP and AP scheme are clearly highlighted: for both Stokes numbers, the AP scheme matches the asymptotic solution whereas the non-AP scheme comes up with a significant deviation which increases while the Stokes number decreases. The quality of the results are also quantitatively assessed in Fig.~\ref{St0.0001_ExplicitError}, in which the error on the density against the number of cells is plotted for the non-AP and the AP schemes and $\St=0.0001$. It shows two orders of magnitude between the two schemes, definitely demonstrating the necessity of AP schemes in such regimes, and the quality of the proposed explicit one.

However, even if the explicit AP scheme obtains good results for every Stokes number, it is constrained by the explicit description of acoustic waves. To avoid such a constraint, the Implicit scheme is a solution. We recall that by Implicit, we mean that the acoustic and source terms are addressed implicitly, where the material wave are still addressed explicitly.

\begin{figure}

\centering
\includegraphics[width=1.0\textwidth]{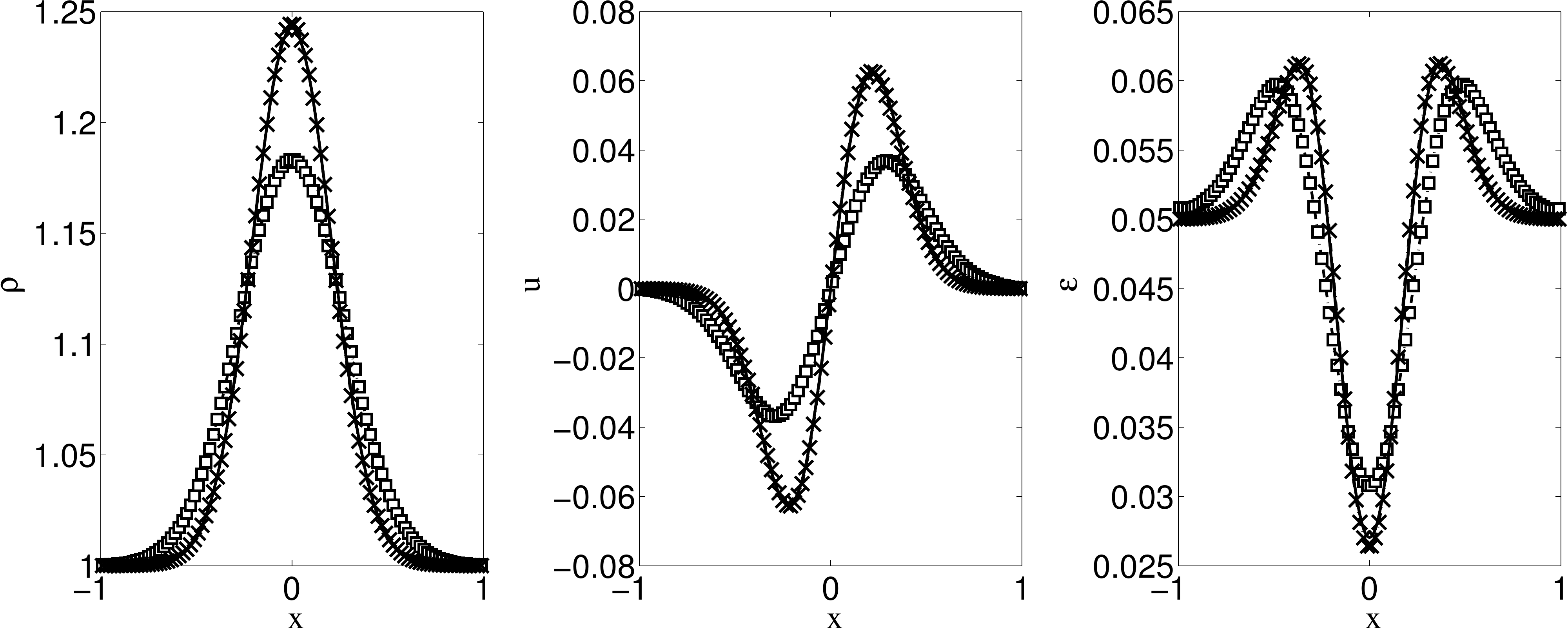}
\caption{Stokes number=0.0001: Comparison between asymptotic analytic solution (black full line), Implicit non-AP scheme (squares) and Implicit AP scheme (crosses) with 100 cells at time $t=0.2s$ and $\Delta t=10\Delta t_{explicit}$.}
\label{St0.001_implicit10dte}
\end{figure}

\begin{figure}

\centering
\includegraphics[width=1.0\textwidth]{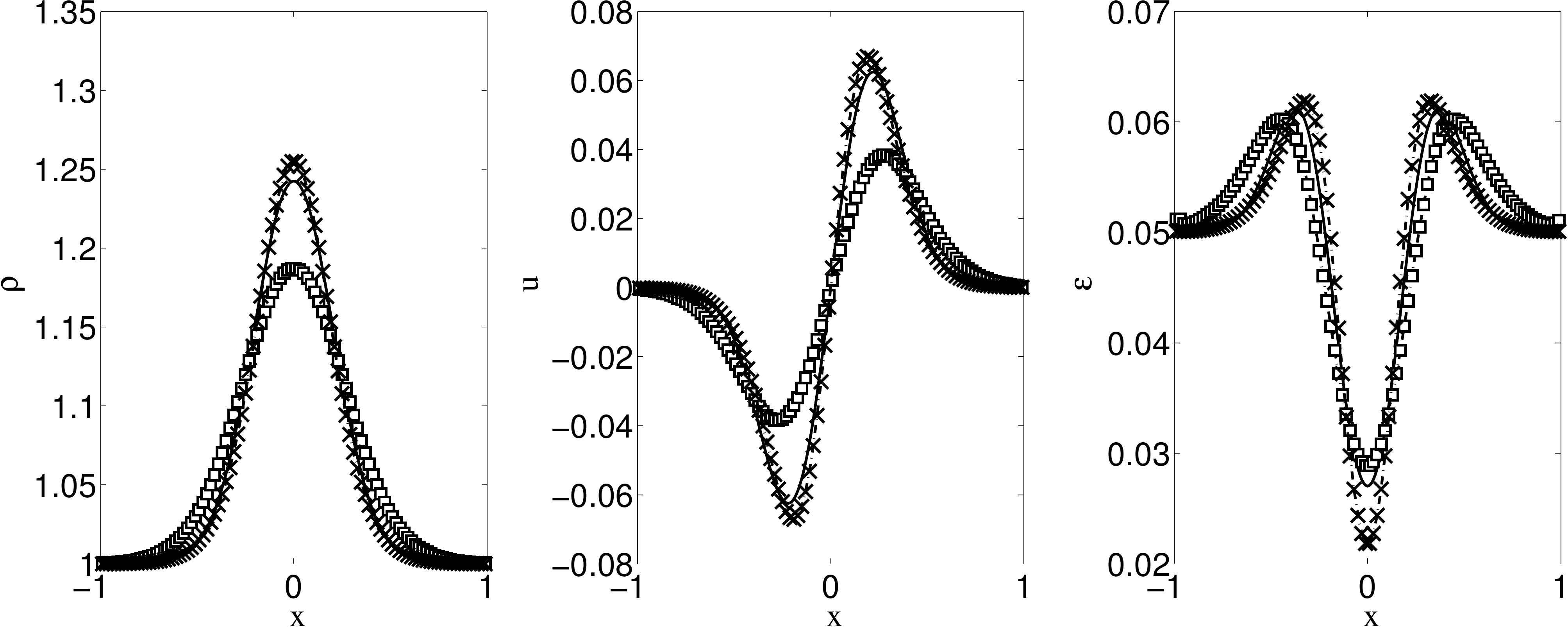}
\caption{Stokes number=0.0001: Comparison between asymptotic analytic solution (black full line), Explicit non-AP scheme (squares) and Implicit AP scheme (crosses) with 100 cells at time $t=0.2s$ and $\Delta t=50\Delta t_{explicit}$.}
\label{St0.001_implicit50dte}
\end{figure}

\begin{figure}

\centering
\includegraphics[width=0.5\textwidth]{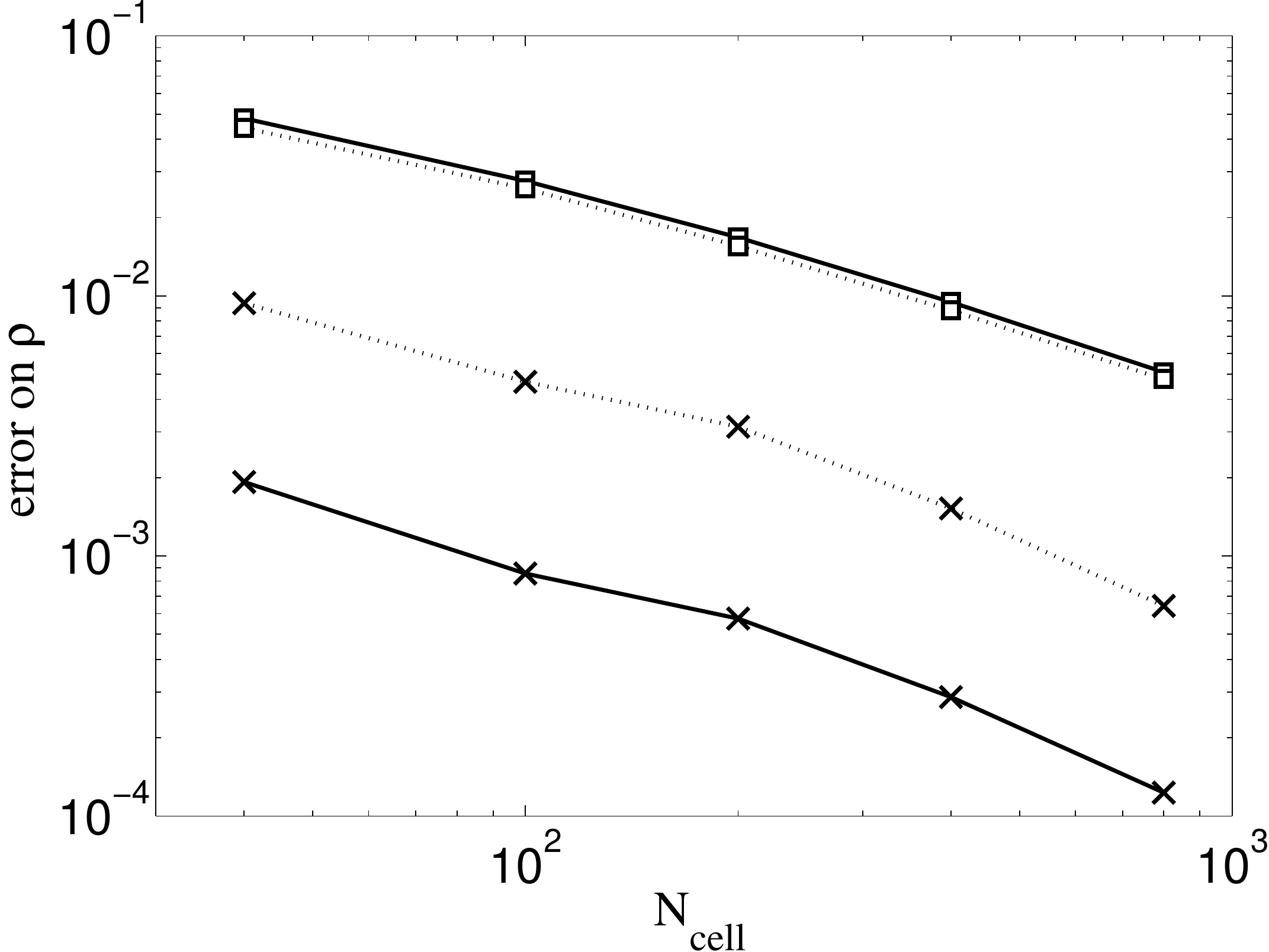}
\caption{Stokes number=0.0001: Error on the number density comparing to the asymptotic analytic solution using the Implicit non-AP scheme (squares) and the Implicit AP scheme (crosses) for $\Delta t=10\Delta t_{explicit}$ (full lines) and $\Delta t=50\Delta t_{explicit}$ (dot lines).}
\label{St0.0001_ImplicitError}
\end{figure}

	In Figs.~\ref{St0.001_implicit10dte}-\ref{St0.001_implicit50dte}, the results for the Implicit AP and Implicit non AP schemes are presented, for $\St=0.0001$ and for a time step which is respectively 10 times and 50 times larger than the explicit time step. It is shown, even using big time steps comparing to the explicit one, increasing the error, the Implicit AP scheme still matches the asymptotic solution, where the non AP one obtains less precise results. The qualitative evaluation of the error on the density plotted in Fig.\ref{St0.0001_ImplicitError} demonstrates it again, the AP scheme being impacted by larger time steps but keeping a low error comparing to the explicit one. 

\FloatBarrier

\section{Conclusions}
In this work, a comparison of existing fully Eulerian strategies for LES has been presented. This comparison highlighted two important aspects of each method: the way the moment equations are closed and the asymptotic behavior at small Stokes number.

Concerning the closure of the moment equations,  on the one side ACBMM can not always ensure the realizability of the moments, because of the lack of a clear link between the moment closure and an underlying NDF. On the other side, KBMM are intrinsically realizable, as long as the underlying presumed NDF is realizable, and accurate and stable numerics are used.


Regarding the available method for particle-laden flow LES, a Kinetic-Based LES method based on the work of \cite{zaichik2009} and \cite{pandya2002} has been chosen. The resulting moment method has two main drawbacks: when the Stokes number tends to zero, the source terms become infinite and in the case of Euler equations, the sound speed too. Moreover the system of equations must tend to a clear advection-diffusion asymptotic limit. To handle it, a new numerical scheme has been designed based on the work of \cite{chalons2012ap} and constituted of a Lagrangian-Projection, a relaxation method and a HLLC scheme with source terms, with a specific treatment of the peculiarities of the models we work with, that is the singular behavior of the fluxes and the inclusion of energy relaxation. The resulting scheme is proven to be Asymptotic-Preserving for the small Stokes number limit in 1D cases. 

The next step is to use the strategy for the simulation of multi-dimensional flows, especially by considering an Anisotropic Gaussian closure \cite{vie2013cicp} to recover the full details of the asymptotic limit. The adaptation of the scheme to space varying Stokes number (for example in evaporating flows) has also to be investigated.

\section{Acknowledgments}
The support of the France-Stanford Center for Interdisciplinary Studies through a collaborative 
project grant  (PIs: P.~Moin and M.~Massot) is also gratefully acknowledged. 
The post-doctoral stay of A.Vi\'e has also been supported by the ANR Sechelles (PIs S. Descombes and M. Massot) and DIGITEO MUSE Project (PI M. Massot).

\bibliographystyle{plain}
\bibliography{biblio}

\end{document}